\documentclass[reqno]{amsart}
\usepackage{tikz,amsmath,amsthm,amsfonts,amssymb,latexsym,mathrsfs,graphicx}

\usepackage{hyperref}

\usepackage{enumerate}
\usepackage[shortlabels]{enumitem}
\usepackage{color}
\usepackage{orcidlink}
\newtheorem{theorem}{Theorem}[section]
\newtheorem{corollary}[theorem]{Corollary}
\newtheorem{lemma}[theorem]{Lemma}

\theoremstyle{definition}
\newtheorem{definition}[theorem]{Definition}

\newtheorem*{problem}{Problem}

\tikzstyle{vertex}=[circle, draw, inner sep=0pt, minimum size=6pt]

\newcommand{\Aut}{{\mathrm {Aut}}}

\newcommand{\id}{\mathrm{id}}

\makeatletter
\newcommand*{\rom}[1]{\expandafter\@slowromancap\romannumeral #1@}
\makeatother
\begin{document}
	
\title{On Subgroup perfect codes and regular sets in generalized Cayley graphs}

\author[Y. Khaefi and F. Sadat Seiedali]
{ Yasamin Khaefi$^{1,*}$ \orcidlink{0009-0006-9468-3093} and Fateme Sadat Seiedali$^{2}$}

\thanks{$^{1}$ Department of Computer Science, Faculty of Mathematics, Physics and
	Informatics, Comenius University, 842 48, Bratislava, Slovakia}

\thanks{$^{*}$ Corresponding author: Yasamin Khaefi.}
\thanks{$^{2}$ Faculty of Mathematics and Computer Science,
	Amirkabir University of Technology (Tehran Polytechnic),
	15914 Tehran, Iran.}

\thanks{E-mail addresses: \href{mailto:yasamin.khaefi@fmph.uniba.sk}{yasamin.khaefi@fmph.uniba.sk} (Y. Khaefi),
	\href{mailto:fseiedali@gmail.com}{fseiedali@gmail.com} (F. S. Seiedali).}
%\MSC[MSC 2010]{05C25 · 05C69 · 94B99}

\begin{abstract}

A subset $C$ of the vertex set of a graph $\Gamma$ is said to be an $(a,b)$-regular set if $C$ induces an $a$-regular subgraph and every vertex outside $C$ is adjacent to exactly $b$ vertices in $C$. A $(0,1)$-regular set is called a perfect code. Let $G$ be a finite group and let $\alpha\in \mathrm{Aut}(G)$ be such that $\alpha^2=\mathrm{id}$. Let
$\omega_\alpha(G)=\{\alpha(g^{-1})g:g\in G\}$.
For a subset $S\subseteq G$ satisfying $\alpha(S)=S^{-1}$ and $S\cap \omega_\alpha(G)=\emptyset$, the generalized Cayley graph of $(G,\alpha)$ with respect to $S$ is the graph with vertex set $G$ in which two distinct elements $g,h\in G$ are adjacent if and only if $\alpha(g^{-1})h\in S$.
In this paper, we prove two main results. First, let $H$ be an $\alpha$-invariant subgroup of a finite group $G$, and let $P$ be an $\alpha$-invariant Sylow $2$-subgroup of $H$. We show that $H$ is a perfect code of $(G,\alpha)$ if and only if $P[H,\alpha]$ is a perfect code of $(G,\alpha)$, where $[H,\alpha]=\langle h^{-1}\alpha(h) : h\in H\rangle$. 
Second, we consider finite groups $G$ possessing a cyclic subgroup $A$ of index $2$ that is invariant under an involutory automorphism $\alpha$. We show that $\omega_\alpha(G)$ is a subgroup of $A$. Using this result, we give a complete characterization of all pairs $(a,b)$ for which an $\alpha$-invariant subgroup $H$ of $G$ is an $(a,b)$-regular set of $(G,\alpha)$. In particular, our characterization applies to all finite groups possessing a characteristic cyclic subgroup of index $2$, for every involutory automorphism $\alpha$.
\end{abstract}

\maketitle
\section{Introduction}
In this paper, all groups are finite and all graphs are simple. Let $\Gamma=({\bf V}(\Gamma),{\bf E}(\Gamma))$ be a graph, where ${\bf V}(\Gamma)$ and ${\bf E}(\Gamma)$ denote its vertex set and edge set, respectively. For non-negative integers $a$ and $b$, a subset $C$ of ${\bf V}(\Gamma)$ is called an $(a,b)$-regular set in $\Gamma$ if every vertex of $C$ has exactly $a$ neighbors in $C$ and every vertex of ${\bf V}(\Gamma)\setminus C$ has exactly $b$ neighbors in $C$; see \cite{R1}. A $(0,1)$-regular set is called a perfect code.

Let $G$ be a finite group, and let $\alpha$ be an automorphism of $G$ such that $\alpha^2=\mathrm{id}$. Define
\[
\begin{aligned}
	T_\alpha(G) &:= \{g\in G : \alpha(g)=g^{-1}\},\\
	\omega_\alpha(G) &:= \{\alpha(g^{-1})g : g\in G\},\\
	\Omega_\alpha(G) &:= T_\alpha(G)\setminus \omega_\alpha(G).
\end{aligned}
\]
Let $S$ be a subset of $G$ such that $S\cap \omega_\alpha(G)=\emptyset$ and $\alpha(S)=S^{-1}$. The generalized Cayley graph of $(G,\alpha)$ with respect to $S$, denoted by $\mathrm{GC}(G,S,\alpha)$, is the graph with vertex set $G$, in which two elements $g,h\in G$ are adjacent if and only if
$
\alpha(g^{-1})h\in S
$.
The condition $S\cap \omega_\alpha(G)=\emptyset$ ensures that $\mathrm{GC}(G,S,\alpha)$ has no loops, while the condition $\alpha(S)=S^{-1}$ ensures that it is undirected. A subset $C$ of $G$ is called a perfect code of $(G,\alpha)$ if there exists a subset $S\subseteq G$ satisfying the above conditions such that $C$ is a perfect code in $\mathrm{GC}(G,S,\alpha)$. Similarly, by an $(a,b)$-regular set of $(G,\alpha)$, we mean an $(a,b)$-regular set in $\mathrm{GC}(G,S,\alpha)$ for some such subset $S$. If a subgroup of $G$ is a perfect code of $(G,\alpha)$, then it is called a subgroup perfect code of $(G,\alpha)$; if it is an $(a,b)$-regular set of $(G,\alpha)$, then it is called a subgroup regular set of $(G,\alpha)$.

When $\alpha=\mathrm{id}$, the conditions $\alpha(S)=S^{-1}$ and $S\cap \omega_\alpha(G)=\emptyset$ reduce to $S=S^{-1}$ and $e\notin S$, respectively. In this case, the generalized Cayley graph $\mathrm{GC}(G,S,\mathrm{id})$ is the usual Cayley graph of $G$ with connection set $S$, and is denoted by $\mathrm{Cay}(G,S)$. Perfect codes in Cayley graphs have been studied in many papers; see, for example, \cite{1, cyclic, 0, Khaefi, Zhang-2}. They are of particular interest, partly because of their connection with Hamming codes \cite{biggs}.

Given the close connection between perfect codes and regular sets, it is natural to study their relationship. In the setting of Cayley graphs, this relationship has been investigated in several papers; see, for example, \cite{Khaefiy, regularset, Wang, Zhang}. For instance, in \cite{Khaefiy}, it is shown that a subgroup $H$ of $G$ is a perfect code of $G$ if and only if it is an $(a,b)$-regular set of $G$ for every pair $(a,b)$ satisfying the required parity conditions.

Perfect codes and regular sets have also been investigated in several graph classes closely related to Cayley graphs, including Cayley sum graphs and coset graphs; see, for example, \cite{Abdol, ma, Seiedalif, wangZhang}. More recently, these notions have also been studied in generalized Cayley graphs; see \cite{Qianfen, seiedali2026generalized}.

In this paper, we pursue this line of research and study perfect codes and regular sets in generalized Cayley graphs. In \cite{seiedali2026generalized}, it is shown that, for a fixed connection set $S$, a subgroup $H$ is a perfect code in $\mathrm{GC}(G,S,\alpha)$ precisely when $\alpha(H)=H$ and $S\cup\{e\}$ forms a transversal of $H$ in $G$. The authors also give a necessary and sufficient condition for a subgroup $H$ of $G$ to be a perfect code of $(G,\alpha)$. 
Building on these results, in Theorem~\ref{L_Theorem}
we show that the subgroup perfect code problem can
be reduced to a certain subgroup $\mathcal{P}$ of $H$.
More precisely, $H$ is a perfect code of $(G,\alpha)$
if and only if $\mathcal{P}$ is a perfect code of
$(G,\alpha)$. 
 This result extends Zhang's reduction theorem for subgroup perfect codes in ordinary Cayley graphs to generalized Cayley graphs; see \cite{Zhang-2}. While Zhang's result reduces the problem to a Sylow $2$-subgroup $P$ of $H$, our approach takes into account the action of $\alpha$ through the subgroup $\mathcal P$. In particular, when $\alpha=\mathrm{id}$, our result recovers Zhang's theorem.

Another result of \cite{seiedali2026generalized} uses the fact that, for an abelian group $G$, the set $\omega_\alpha(G)$ is a subgroup of $G$, and characterizes all pairs $(a,b)$ for which a subgroup $H$ of $G$ is an $(a,b)$-regular set of $(G,\alpha)$. In the final section of this paper, we extend this study to finite groups possessing a cyclic subgroup $A$ of index $2$ that is invariant under an involutory automorphism $\alpha$. Unlike the abelian case, in our result the subgroup $H$ need not be normal in $G$. By proving that $\omega_\alpha(G)$ is a subgroup and handling the resulting coset structure, we give a complete characterization of the possible pairs $(a,b)$; see Lemma~\ref{omega_sub_new} and Theorem~\ref{thm:uni-new}. In particular, our characterization applies to all finite groups with a characteristic cyclic subgroup of index $2$, for every involutory automorphism $\alpha$. These include dihedral, dicyclic, and semidihedral groups, among others. Moreover, our result also applies to modular groups, provided that their standard cyclic subgroup of index $2$ is $\alpha$-invariant.

\section{Preliminaries}\label{sec:prelim}
We begin with two preliminary results from \cite{seiedali2026generalized}.
\begin{lemma}\cite[Lemma~2.7]{seiedali2026generalized}\label{2.7}
	Subgroup $H$ is an $(a,b)$-regular set in $GC(G,S,\alpha)$ if and only if $S=S_1\cup S_2$ such that
	$S_1\subseteq H$, $|S_1|=a$, $S_2$ is a union of $b$ pairwise disjoint right transversals of $H$ in $G\setminus H$,
	and $H$ is an $(a,0)$-regular set in $GC(G,S_1,\alpha)$ and a $(0,b)$-regular set in $GC(G,S_2,\alpha)$.
\end{lemma}
\begin{lemma}\cite[Theorem~A]{seiedali2026generalized}\label{ThA}
	Let $G$ be a group and $H$ a subgroup of $G$. Let $\alpha \in \operatorname{Aut}(G)$ with $\alpha^2 = \operatorname{id}$ such that $\alpha(H) = H$. Then $H$ is a subgroup perfect code of $(G,\alpha)$ if and only if there is no $g \in G \setminus H$ such that 
	$
	\alpha(HgH) = Hg^{-1}H$, $|HgH|/|H|$  is odd, and $HgH \cap \Omega_{\alpha}(G) = \emptyset$.
\end{lemma}
\section{Subgroup perfect codes}

Let $H\leq G$ be an $\alpha$-invariant subgroup, where
$\alpha\in\operatorname{Aut}(G)$ and $\alpha^2=\id$.
Since the number of Sylow $2$-subgroups of $H$ is odd
and $\alpha$ is an involutory automorphism, there exists
an $\alpha$-invariant Sylow $2$-subgroup $P$ of $H$.
We fix such a subgroup $P$ and define
\[
[H,\alpha]
=\langle h^{-1}\alpha(h):h\in H\rangle,
\qquad
\mathcal{P}=P[H,\alpha].
\]

Although $\mathcal{P}$ may depend on the choice of $P$, our results hold for every $\alpha$-invariant Sylow $2$-subgroup of $H$. Thus, we fix one such subgroup throughout and refer to its associated $\mathcal{P}$ as the {\em $\alpha$-root} of $H$. This definition is motivated by the ordinary Cayley
graph case ($\alpha=\id$). In this case,
$[H,\alpha]=1$, and hence $\mathcal{P}=P$ is a
Sylow $2$-subgroup of $H$. Recall that, for ordinary
Cayley graphs, a subgroup $H$ of $G$ is a perfect code
of some Cayley graph of $G$ if and only if a Sylow
$2$-subgroup of $H$ is a perfect code of some Cayley
graph of $G$; see \cite[Theorem 1.1]{Zhang-2}.
Our aim is to extend this reduction to generalized
Cayley graphs by taking into account the action of
$\alpha$ through the subgroup $[H,\alpha]$.

In the following, we first recall two group-theoretic lemmas and establish some basic properties of \(\mathcal P\) needed for our main result.

\begin{lemma}\cite[Glauberman 3.24]{Isaacs}\label{Glauberman}
	Let $A$ act via automorphisms on $G$, where  $A$
	and  $G$ are finite groups, and  suppose that
	$
	(|A|,|G|)=1
	$.
	Assume also that at least one of $A$ or $G$ is solvable. Suppose that 
	$A$ and $G$ each act on some nonempty set $X$, where the action of $G$
	 is transitive. Finally, assume the
	compatibility condition:
	\[
	(\alpha\cdot g)\cdot a
	=
	(\alpha\cdot a)\cdot g^a
	\]
	for all $\alpha\in X$, $a\in A$, and $g\in G$. Then:
	\begin{enumerate}
		\item[(a)] There exists an  $A$-invariant element $\alpha\in X$.
		\item[(b)] If $\alpha,\beta\in X$ are $A$-invariant, then
		there exists $c\in C_G(A)$ such that
		$
		\alpha\cdot c=\beta
		$.
	\end{enumerate}
\end{lemma}

\begin{lemma}\label{basic}
    Let \(G\) be a  group, let \(H\le G\), and let \(P\) be a Sylow \(2\)-subgroup of \(H\). Then the following hold:
    \begin{enumerate}
        \item If 
$[H:H\cap xHx^{-1}]$ is odd for some $x\in G$, then there exists \(y\in HxH\) such that \(y\in \mathrm{N}_G(P)\).
\item  If \(y\in \mathrm{N}_G(P)\), then $HyH\cap \mathrm{N}_G(P)=\mathrm{N}_H(P)y\mathrm{N}_H(P)$.
    \end{enumerate}
\end{lemma}

\begin{proof}
(1) Put $I=H\cap xHx^{-1}$ and let
$Q\in\operatorname{Syl}_2(I)$. Since
$[H:I]$ is odd, $Q$ is a Sylow
$2$-subgroup of  $H$ and $xHx^{-1}$.
Thus there exist $h,k\in H$ such that
$Q=hPh^{-1}=xkPk^{-1}x^{-1}$.
Hence $h^{-1}xk\in HxH\cap \mathrm{N}_G(P)$.

(2) Clearly, $\mathrm{N}_H(P)y\mathrm{N}_H(P)\subseteq HyH\cap \mathrm{N}_G(P)$. Let $y\in \mathrm{N}_G(P)$ and take
$z=h_1yh_2\in HyH\cap N_G(P)$,
where $h_1,h_2\in H$. Then $P$ and
$h_1^{-1}Ph_1$ are Sylow $2$-subgroups
of $I=H\cap yHy^{-1}$.
Hence there exists $u\in I$ such that
$u^{-1}h_1^{-1}Ph_1u=P$.
Thus $h_1u\in \mathrm{N}_H(P)$, and
$z=(h_1u)y(y^{-1}u^{-1}yh_2)
\in \mathrm{N}_H(P)y\mathrm{N}_H(P)$.
\end{proof}

\begin{lemma}\label{H,alpha}
    Let $G$ be a group and $\alpha$ be an involutory automorphism of $G$. Let $H$ be an $\alpha$-invariant subgroup of $G$. Then $[H,\alpha]\trianglelefteq H$, $\alpha([H,\alpha])=[H,\alpha]$, and $\alpha$ induces the identity automorphism on $H/[H,\alpha]$.
\end{lemma}
\begin{proof}
    Let $n_h=h^{-1}\alpha(h)$, for each $h\in H$. Then $x^{-1}n_hx=(hx)^{-1}\alpha(hx)(x^{-1}\alpha(x))^{-1}\in [H,\alpha]$ for every $x\in H$, and so $[H,\alpha]$ is a normal subgroup of $H$. Clearly $\alpha(n_h)=n_h^{-1}$ and so $\alpha([H,\alpha])=[H,\alpha]$. Since $h^{-1}\alpha(h)\in [H,\alpha]$, for every $h\in H$, we have $\alpha(h)[H,\alpha]=h[H,\alpha]$ and therefore $\alpha$ is the identity automorphism of $H/[H,\alpha]$.
\end{proof}

\begin{lemma}\label{mathcal(P)}
	Let $G$ be a group and $\alpha$ be an involutory automorphism of $G$. Let $H$ be an $\alpha$-invariant subgroup of $G$ and $\mathcal{P}$ be the $\alpha$-root of $H$. Then the following hold:
    \begin{enumerate}
        \item $\mathcal{P}$ is an $\alpha$-invariant subgroup of $H$ and $[H:\mathcal{P}]$ is odd; 
        \item  $\mathrm{N}_{\mathcal{P}}(P)\trianglelefteq \mathrm{N}_H(P)$, $[\mathrm{N}_H(P): \mathrm{N}_{\mathcal{P}}(P)]$ is odd, and $\alpha$ induces the  identity automorphism on $\mathrm{N}_H(P)/\mathrm{N}_{\mathcal{P}}(P)$.
    \end{enumerate}
\end{lemma}

\begin{proof}
(1) By Lemma~\ref{H,alpha},
$[H,\alpha]\trianglelefteq H$ and
$\alpha([H,\alpha])=[H,\alpha]$.
Hence
$\mathcal P=P[H,\alpha]
$
is an $\alpha$-invariant subgroup
of $H$. Since $P$ is a Sylow
$2$-subgroup of $H$ and
$P\leq\mathcal P\leq H$,
the index $[H:\mathcal P]$ is odd.

(2) Put $N=[H,\alpha]$ and
$B=\mathrm{N}_H(P)$. Since $P\leq B$,
$
\mathrm{N}_{\mathcal P}(P)
=\mathcal P\cap B
=PN\cap B
=P(N\cap B)$.
Since $P\trianglelefteq B$ and
$N\cap B\trianglelefteq B$, we get that
$\mathrm{N}_{\mathcal P}(P)\trianglelefteq B.$
Moreover, $P$ is a Sylow
$2$-subgroup of $B$, and hence
$[B:\mathrm{N}_{\mathcal P}(P)]
$
is odd. Finally, $\alpha(B)=B$,
and for every $x\in B$,
$x^{-1}\alpha(x)\in N\cap B
\leq \mathrm{N}_{\mathcal P}(P).
$
Thus $\alpha$ induces the identity
automorphism on
$B/\mathrm{N}_{\mathcal P}(P)$.
\end{proof}

%%%%%%%%%%%%%%%%%%%
%%%%%%%%%%%%%%%%%%%

\begin{lemma}\label{Omega}
	Let $G$ be a group and $\alpha$ be an involutory automorphism of $G$. 
	Let $H$ be an $\alpha$-invariant subgroup of
	$G$ and  $\mathcal{P}$ be the $\alpha$-root of $H$. Suppose that
	$x\in G\setminus H$ satisfies $\alpha(\mathcal{P}x\mathcal{P})=\mathcal{P}x^{-1}\mathcal{P}$ and $|\mathcal{P}x\mathcal{P}|/|\mathcal{P}|$ is odd.
	 If $HxH\cap \Omega_\alpha(G)\ne\emptyset$, then  $\mathcal{P}x\mathcal{P}\cap \Omega_\alpha(G)\ne\emptyset$.
\end{lemma}

\begin{proof}
As $|\mathcal Px\mathcal P|/|\mathcal P|$
is odd, by Lemma \ref{basic},
there exists
$
y\in\mathcal{P}x\mathcal{P}
\cap \mathrm{N}_G(P)$.
Thus
$\mathcal{P}y\mathcal{P}
=\mathcal{P}x\mathcal{P}$, and 
$HyH=HxH$.
By Lemma~\ref{mathcal(P)}, we have
$\mathrm{N}_{\mathcal{P}}(P)\trianglelefteq \mathrm{N}_H(P)$, and
$[\mathrm{N}_H(P):\mathrm{N}_{\mathcal{P}}(P)]$ is odd,
and both $\mathrm{N}_{\mathcal{P}}(P)$ and $\mathrm{N}_H(P)$ are
$\alpha$-invariant.
Since
$\alpha(\mathcal{P}y\mathcal{P})
=\mathcal{P}y^{-1}\mathcal{P}$,
we have
$\alpha(y)^{-1}
\in\mathcal{P}y\mathcal{P}$.
Moreover, $\alpha(P)=P$ and
$y\in \mathrm{N}_G(P)$ imply that
$\alpha(y)^{-1}\in \mathrm{N}_G(P)$.
Hence, by
Lemma~\ref{basic},
$\alpha(y)^{-1}\in \mathrm{N}_{\mathcal{P}}(P)y\mathrm{N}_{\mathcal{P}}(P)$.
Therefore,
$\alpha(\mathrm{N}_{\mathcal{P}}(P)y\mathrm{N}_{\mathcal{P}}(P))=\mathrm{N}_{\mathcal{P}}(P)y^{-1}\mathrm{N}_{\mathcal{P}}(P)$.
Now we choose
$n\in HyH\cap\Omega_\alpha(G)$ and set
$I=H\cap nHn^{-1}$.
Since $n\in HyH$ and
$y\in \mathrm{N}_G(P)$, the index
$[H:I]$ is odd. In particular,
every Sylow $2$-subgroup of $I$
is also a Sylow $2$-subgroup of $H$.
We define a map
$\theta:I\longrightarrow I$, as
$\theta(h)=n\alpha(h)n^{-1}$.
Since $\alpha(n)=n^{-1}$,
 $\theta$ is an
involutory automorphism of $I$.
As the number of Sylow
$2$-subgroups of $I$ is odd,
there exists a
$\theta$-invariant Sylow
$2$-subgroup $Q$ of $I$.
Since $Q$ and $P$ are Sylow
$2$-subgroups of $H$, there
exists $a\in H$ such that
$Q=aPa^{-1}$.
Put
$m=a^{-1}n\alpha(a)$.
Using $\theta(Q)=Q$ and
$\alpha(P)=P$, we get that
$mPm^{-1}=P$.
Thus $m\in N_G(P)$.
Moreover, 
$a^{-1}\Omega_\alpha(G)\alpha(a)=\Omega_\alpha(G)$, for each $a\in G$. So
 $m\in\Omega_\alpha(G)$.
Since $m\in HyH$, we get that
$m\in HyH\cap N_G(P)
\cap\Omega_\alpha(G)$.
Then by Lemma~\ref{basic},
$m\in \mathrm{N}_H(P)y\mathrm{N}_H(P)\cap\Omega_\alpha(G)$.
Now we apply Lemma~\ref{Glauberman}.
Let
\[
\mathcal{S}
=\{\mathrm{N}_{\mathcal{P}}(P)z\mathrm{N}_{\mathcal{P}}(P):z\in \mathrm{N}_H(P)y\mathrm{N}_H(P)\}
\]
and
\[
D=(\mathrm{N}_H(P)/\mathrm{N}_{\mathcal{P}}(P))\times(\mathrm{N}_H(P)/\mathrm{N}_{\mathcal{P}}(P)).
\]
The group $D$ has odd order
and acts transitively on
$\mathcal{S}$ by
\[
(u\mathrm{N}_{\mathcal{P}}(P),v\mathrm{N}_{\mathcal{P}}(P))\cdot \mathrm{N}_{\mathcal{P}}(P)z\mathrm{N}_{\mathcal{P}}(P)
=\mathrm{N}_{\mathcal{P}}(P)uzv^{-1}\mathrm{N}_{\mathcal{P}}(P).
\]
We define an automorphism
$\varphi$ of $D$ by
\[
\varphi(u\mathrm{N}_{\mathcal{P}}(P),v\mathrm{N}_{\mathcal{P}}(P))
=(\alpha(v)\mathrm{N}_{\mathcal{P}}(P),\alpha(u)\mathrm{N}_{\mathcal{P}}(P)).
\]
Let $\varphi$ act on
$\mathcal{S}$ by
\[
\varphi\cdot \mathrm{N}_{\mathcal{P}}(P)z\mathrm{N}_{\mathcal{P}}(P)
=\mathrm{N}_{\mathcal{P}}(P)\alpha(z)^{-1}\mathrm{N}_{\mathcal{P}}(P).
\]
These actions satisfy the
compatibility condition
of Lemma~\ref{Glauberman}.
Also both $\mathrm{N}_{\mathcal{P}}(P)y\mathrm{N}_{\mathcal{P}}(P)$ and $\mathrm{N}_{\mathcal{P}}(P)m\mathrm{N}_{\mathcal{P}}(P)$ are
fixed by $\varphi$. Indeed,
$\alpha(\mathrm{N}_{\mathcal{P}}(P)y\mathrm{N}_{\mathcal{P}}(P))=\mathrm{N}_{\mathcal{P}}(P)y^{-1}\mathrm{N}_{\mathcal{P}}(P)$
and $\alpha(m)=m^{-1}$.
Since $|D|$ is odd,
Lemma~\ref{Glauberman} implies
that $C_D(\varphi)$ acts
transitively on the
$\varphi$-fixed elements
of $\mathcal{S}$.
Therefore there exists
$(u\mathrm{N}_{\mathcal{P}}(P),v\mathrm{N}_{\mathcal{P}}(P))\in C_D(\varphi)$
such that
$\mathrm{N}_{\mathcal{P}}(P)m\mathrm{N}_{\mathcal{P}}(P)=\mathrm{N}_{\mathcal{P}}(P)uyv^{-1}\mathrm{N}_{\mathcal{P}}(P)$.
Since
$(u\mathrm{N}_{\mathcal{P}}(P),v\mathrm{N}_{\mathcal{P}}(P))\in C_D(\varphi)$,
we have
$v\mathrm{N}_{\mathcal{P}}(P)=\alpha(u)\mathrm{N}_{\mathcal{P}}(P)$ and so
$\mathrm{N}_{\mathcal{P}}(P)m\mathrm{N}_{\mathcal{P}}(P)=\mathrm{N}_{\mathcal{P}}(P)uy\alpha(u)^{-1}\mathrm{N}_{\mathcal{P}}(P)$.
As $\mathrm{N}_{\mathcal{P}}(P)\trianglelefteq \mathrm{N}_H(P)$,
it follows that
$u^{-1}m\alpha(u)\in \mathrm{N}_{\mathcal{P}}(P)y\mathrm{N}_{\mathcal{P}}(P)$.
Moreover,
$m\in\Omega_\alpha(G)$,
and $\Omega_\alpha(G)$ is
invariant under twisted
conjugation. Hence
$u^{-1}m\alpha(u)
\in \mathrm{N}_{\mathcal{P}}(P)y\mathrm{N}_{\mathcal{P}}(P)\cap\Omega_\alpha(G)$.
Since
$\mathrm{N}_{\mathcal{P}}(P)y\mathrm{N}_{\mathcal{P}}(P)\subseteq
\mathcal{P}y\mathcal{P}
=\mathcal{P}x\mathcal{P}$,
we conclude that
$\mathcal{P}x\mathcal{P}
\cap\Omega_\alpha(G)
\neq\emptyset$.
\end{proof}

We are now ready to prove the main theorem.

\begin{theorem}\label{L_Theorem}
	Let $G$ be a group and  $\alpha$ be an involutory automorphism of
	$G$. Let $H$ be an $\alpha$-invariant subgroup of $G$ and  $\mathcal{P}$ be
	the $\alpha$-root of $H$. Then $H$ is a perfect code of
	$(G,\alpha)$ if and only if $\mathcal{P}$ is a perfect code of $(G,\alpha)$.
\end{theorem}

\begin{proof}
First suppose that $\mathcal{P}$ is a
perfect code of $(G,\alpha)$.
Assume, to the contrary, that $H$ is
not a perfect code of $(G,\alpha)$.
By Lemma~\ref{ThA}, there exists
$x\in G\setminus H$ such that
$\alpha(HxH)=Hx^{-1}H$,
$[H:H\cap xHx^{-1}]$ is odd, and
$HxH\cap\Omega_\alpha(G)=\emptyset$.
Let
$\mathcal{A}
=\{\mathcal{P}y\mathcal{P}:y\in HxH\}$.
Since $\alpha(\mathcal{P})=\mathcal{P}$,
the map
$\tau:\mathcal{A}\longrightarrow\mathcal{A}$,
with
$\tau(T)=\alpha(T^{-1})$, for $T\in \mathcal{A}$ is a well-defined involution.
By Lemma~\ref{mathcal(P)},
$[H:\mathcal{P}]$ is odd. Hence
$
\sum_{T\in\mathcal{A}}
|T|/|\mathcal{P}|=
|HxH|/|\mathcal{P}|
=[H:\mathcal{P}]
  [H:H\cap xHx^{-1}]$
is odd.
Since $\tau$ preserves the size of
each double coset, there exists a
$\tau$-fixed double coset
$\mathcal{P}y\mathcal{P}$
such that
$|\mathcal{P}y\mathcal{P}|/|\mathcal{P}|
=[\mathcal{P}:
\mathcal{P}\cap y\mathcal{P}y^{-1}]
$
is odd.
As $\tau(\mathcal{P}y\mathcal{P})=\mathcal{P}y\mathcal{P}$, we have
$\alpha(\mathcal{P}y\mathcal{P})
=\mathcal{P}y^{-1}\mathcal{P}$.
Moreover, $y\notin\mathcal{P}$,
because $y\in HxH$ and $x\notin H$.
Since $\mathcal{P}y\mathcal{P}\subseteq HxH$, we also have
$\mathcal{P}y\mathcal{P}\cap\Omega_\alpha(G)=\emptyset$.
This contradicts Lemma~\ref{ThA},
since $\mathcal{P}$ is a perfect
code of $(G,\alpha)$.
Therefore $H$ is a perfect code
of $(G,\alpha)$.

Conversely, suppose that $H$ is
a perfect code of $(G,\alpha)$.
Assume, to the contrary, that
$\mathcal{P}$ is not a perfect code
of $(G,\alpha)$.
By Lemma~\ref{ThA}, there exists
$x\in G\setminus\mathcal{P}$ such that
$\alpha(\mathcal{P}x\mathcal{P})
=\mathcal{P}x^{-1}\mathcal{P}$,
$[\mathcal{P}:
\mathcal{P}\cap x\mathcal{P}x^{-1}]$
is odd, and
$\mathcal{P}x\mathcal{P}
\cap\Omega_\alpha(G)=\emptyset$.
By Lemma~\ref{basic},
there exists
$y\in\mathcal{P}x\mathcal{P}
\cap\mathrm{N}_G(P)$. Thus
$\mathcal{P}y\mathcal{P}
=\mathcal{P}x\mathcal{P}$.
By Lemma~\ref{mathcal(P)}, we have
$\mathrm{N}_{\mathcal{P}}(P)\trianglelefteq \mathrm{N}_H(P)$, $[\mathrm{N}_H(P):\mathrm{N}_{\mathcal{P}}(P)]$ is odd,
and
 $\alpha$ induces the identity
automorphism on $\mathrm{N}_H(P)/\mathrm{N}_{\mathcal{P}}(P)$.
Since
$\alpha(\mathcal{P}y\mathcal{P})
=\mathcal{P}y^{-1}\mathcal{P}$,
we have
$\alpha(y)^{-1}
\in\mathcal{P}y\mathcal{P}$.
Moreover, since $\alpha(P)=P$ and
$y\in\mathrm{N}_G(P)$, we get that
$\alpha(y)^{-1}\in\mathrm{N}_G(P)$.
Thus, by
Lemma~\ref{basic},
$\alpha(y)^{-1}\in \mathrm{N}_{\mathcal{P}}(P)y\mathrm{N}_{\mathcal{P}}(P)$.
Therefore,
$\alpha(\mathrm{N}_{\mathcal{P}}(P)y\mathrm{N}_{\mathcal{P}}(P))=\mathrm{N}_{\mathcal{P}}(P)y^{-1}\mathrm{N}_{\mathcal{P}}(P)$.
We claim that $y\notin H$.
Suppose otherwise. Since
$y\in\mathrm{N}_G(P)$, we have
$y\in \mathrm{N}_H(P)$.
As $\mathrm{N}_{\mathcal{P}}(P)\trianglelefteq \mathrm{N}_H(P)$, we get that
$\alpha(y)\mathrm{N}_{\mathcal{P}}(P)=y^{-1}\mathrm{N}_{\mathcal{P}}(P)$.
On the other hand, $\alpha$ acts
trivially on $\mathrm{N}_H(P)/\mathrm{N}_{\mathcal{P}}(P)$, so
$\alpha(y)\mathrm{N}_{\mathcal{P}}(P)=y\mathrm{N}_{\mathcal{P}}(P)$.
Therefore
$(y\mathrm{N}_{\mathcal{P}}(P))^2=\mathrm{N}_{\mathcal{P}}(P)$.
Since $|\mathrm{N}_H(P)/\mathrm{N}_{\mathcal{P}}(P)|$ is odd, it follows
that $y\in \mathrm{N}_{\mathcal{P}}(P)\leq\mathcal{P}$, a contradiction.
Hence $y\notin H$.
Since $y\in\mathrm{N}_G(P)$, we have
$P\leq H\cap yHy^{-1}$.
As $P$ is a Sylow $2$-subgroup
of $H$, it follows that
$[H:H\cap yHy^{-1}]$ is odd.
Moreover,
$\alpha(HyH)=Hy^{-1}H$.
Since $H$ is a perfect code
of $(G,\alpha)$, Lemma~\ref{ThA}
implies that
$HyH\cap\Omega_\alpha(G)
\neq\emptyset$.
Using Lemma~\ref{Omega},
we get that
$\mathcal{P}y\mathcal{P}
\cap\Omega_\alpha(G)
\neq\emptyset$.
This is a contradiction, since
$\mathcal{P}y\mathcal{P}
=\mathcal{P}x\mathcal{P}$. 
Therefore $\mathcal{P}$ is a
perfect code of $(G,\alpha)$.
\end{proof}

\begin{corollary}
Let $G$ be a group, let
$\alpha\in\operatorname{Aut}(G)$ be an involutory automorphism of $G$, and let $H$ be an
$\alpha$-invariant subgroup of $G$.
If $P_1$ and $P_2$ are $\alpha$-invariant
Sylow $2$-subgroups of $H$, then
$P_1[H,\alpha]$ is a perfect code of
$(G,\alpha)$ if and only if
$P_2[H,\alpha]$ is a perfect code of
$(G,\alpha)$.
\end{corollary}

\begin{proof}
By Theorem~\ref{L_Theorem}, both
conditions are equivalent to $H$ being
a perfect code of $(G,\alpha)$.
\end{proof}

In particular, when $\alpha=\id$, we have
$[H,\alpha]=1$, and every Sylow
$2$-subgroup of $H$ is
$\alpha$-invariant. Consequently,
our theorem recovers Zhang's
characterization \cite[Theorem 1.1]{Zhang-2}:
$H$ is a perfect code of $G$ if and
only if some (equivalently, every)
Sylow $2$-subgroup of $H$ is a
perfect code of $G$.

%%%%%%%%%%%%%%%%%%%%%%%%%%%%%%%%%%%
\section{Some regular sets}

The regularity of the subgroups of abelian groups was classified in \cite{seiedali2026generalized}. We now consider finite groups possessing a cyclic subgroup of index $2$. Let $G$ be such a group and let $A=\langle r\rangle$ be a cyclic subgroup of $G$ of index $2$. We show that, for every involutory automorphism $\alpha\in\operatorname{Aut}(G)$ satisfying $\alpha(A)=A$, the set $\omega_\alpha(G)$ is a subgroup of $A$; see Lemma~\ref{omega_sub_new}. Using this result, we determine all possible pairs $(a,b)$ for which an $\alpha$-invariant subgroup $H$ of $G$ is an $(a,b)$-regular set of $(G,\alpha)$.

%%%%%

In the following, for a proper subgroup $H$ of $G$ and
each $x\in G\setminus H$, we set
$\mathcal{L}_x=\bigl|Hx\setminus \omega_\alpha(G)\bigr|$, and define $ \mathcal{L}_{\min}(H)=\min\{\mathcal{L}_x \mid x\in G\setminus H\}$.

The next lemma is a slight reformulation of \cite[ Lemma~2.9]{seiedali2026generalized}. Although it is stated there for abelian groups, its proof uses only the fact that $\omega_\alpha(G)$ is a subgroup of $G$. Hence the same argument gives the following statement. We include the details for completeness.

\begin{lemma}\label{L(x)}
	Let $G$ be a group and 
$\alpha\in\operatorname{Aut}(G)$ such that
$\alpha^2=\id$. Suppose that
$\omega_\alpha(G)$ is a subgroup of $G$.
Let $H$ be a proper subgroup of $G$. Then:
	\begin{enumerate}[label=(\alph*)]
		\item $\mathcal{L}_{\min}(H)=|H|$ if $\omega_\alpha(G)\subseteq H$;
		\item $\mathcal{L}_{\min}(H)=\bigl|H\setminus \omega_\alpha(G)\bigr|$ if $\omega_\alpha(G)\nsubseteq H$.
	\end{enumerate}
\end{lemma}

	\begin{proof}
		For each $x\in G\setminus H$, we have $\mathcal{L}_x=|Hx\setminus\omega_\alpha(G)|=|Hx|-|Hx\cap\omega_\alpha(G)|$. Suppose that $Hx\cap\omega_\alpha(G)\neq\emptyset$ and take $hx\in Hx\cap\omega_\alpha(G)$, for some $h\in H$. Without loss of generality, we assume that $hx=x$, as $\mathcal{L}_x=\mathcal{L}_{hx}$. Then $|Hx\cap\omega_\alpha(G)|=|H\cap\omega_\alpha(G)x^{-1}|=|H\cap\omega_\alpha(G)|$. Thus either $\mathcal{L}_x=|H|$ or $\mathcal{L}_x=|H|-|H\cap\omega_\alpha(G)|=|H\setminus\omega_\alpha(G)|$.
		
		If $\omega_\alpha(G)\subseteq H$, then $\mathcal{L}_x=|H|$ for all $x\in G\setminus H$, and hence $\mathcal{L}_{\min}(H)=\min\{\mathcal{L}_x:x\in G\setminus H\}=|H|$.
		
		If $\omega_\alpha(G)\nsubseteq H$, then there exists $x\in G\setminus H$ such that $x\in\omega_\alpha(G)$ and so, according to the above argument, $\mathcal{L}_x=|H\setminus\omega_\alpha(G)|\leq |H|$. Hence $\mathcal{L}_{\min}(H)=\min\{\mathcal{L}_x:x\in G\setminus H\}=|H\setminus\omega_\alpha(G)|$.
	\end{proof}

In \cite[Lemma~2.10]{seiedali2026generalized}, the authors also prove the following for abelian groups, though in fact the statement holds for every group. We include the details for completeness.

\begin{lemma}\label{a-regularity}
Let $G$ be a group and
$\alpha\in\operatorname{Aut}(G)$ such that
$\alpha^2=\id$. Let $H$ be an
$\alpha$-invariant subgroup of $G$.
For every integer
$0\leq a\leq |H\setminus\omega_\alpha(G)|$,
the subgroup $H$ is an $(a,0)$-regular set
of $(G,\alpha)$ if and only if
at least one of the following conditions holds:
\begin{enumerate}
\item $a$ is even;
\item $H\cap\Omega_\alpha(G)\neq\emptyset$.
\end{enumerate}
\end{lemma}

\begin{proof}
Put $
W=H\setminus T_\alpha(G)$ and
$F=H\cap\Omega_\alpha(G)$.
Let $\iota(g)=\alpha(g^{-1})$.
Since
$\omega_\alpha(G)\subseteq T_\alpha(G)$,
we have
$H\setminus\omega_\alpha(G)
=W\;\dot\cup\;F$.
The involution $\iota$ fixes
every element of $F$ and acts
without fixed points on $W$.
Hence $W$ is a disjoint union
of two-element orbits.

Suppose that $H$ is an
$(a,0)$-regular set of $(G,\alpha)$.
By Lemma~\ref{2.7}, there exists
an $\iota$-invariant subset
$S\subseteq W\cup F$ of size $a$.
If $F=\emptyset$, then $S$ is
a union of two-element orbits,
so $a$ is even.

Conversely, suppose that $a$
is even or $F\neq\emptyset$.
If $a\leq |W|$ and $a$ is even,
choose $a/2$ two-element orbits
of $W$. If $a\leq |W|$ and $a$
is odd, choose one element of
$F$ together with $(a-1)/2$
two-element orbits of $W$.
If $a>|W|$, choose all elements
of $W$ together with any
$a-|W|$ elements of $F$.
In every case, the resulting
set $S$ is $\iota$-invariant,
is contained in
$H\setminus\omega_\alpha(G)$,
and has cardinality $a$.
Thus $H$ is an $(a,0)$-regular
set by Lemma~\ref{2.7}.
\end{proof}

\begin{lemma}\label{omega_sub_new}
    Let $G$ be a group and $A$ be a cyclic subgroup of $G$ of index 2, and let $\alpha\in\operatorname{Aut}(G)$ be an involutory automorphism such that $\alpha(A)=A$. Then $\omega_\alpha(G)\le A$.
\end{lemma}

\begin{proof}
    Let $A=\langle r \rangle$ and $|r|=n$. Since $[G:A] =2$, we get that $A$ is a normal subgroup of $G$ and there exists $s\in G\setminus A$ such that $G=A\cup As$.  Therefore every element of $G$ is either $r^i$ or $r^is$, where
	$0\le i<n$. As $A$ is normal, we have $s^{-1}rs=r^t$, where $\gcd(t,n)=1$, and $\alpha(A)=A$ implies that $\alpha(r)=r^k$, where $\gcd(k,n)=1$. Also $\alpha(As)=As$ and so $ \alpha(s)=r^m s,
	$ for some integer $m$. Moreover $\alpha^2(s)=s$ implies that $s=r^{m(k+1)}s$ and so $m(k+1)\equiv 0 \pmod{n}$.
    Let $d=\gcd(n,1-k)$. Now we have two cases:
\begin{itemize}
		\item For $g=r^i$,
		$\alpha(g^{-1})g=
		r^{(1-k)i}$ and so 
        $\{r^{(1-k)i}\mid 0\le i<n\}=\langle r^{1-k}\rangle = \langle r^d\rangle$.

		\item For $g=r^is$, $\alpha(g^{-1})g= 
        r^{t(-m+(1-k)i)}$. 
         We have $\langle r^{t(1-k)}\rangle= \langle r^{1-k}\rangle=\langle r^d\rangle$, since $\gcd (t,n)=1$. So 
         $\{r^{-tm+t(1-k)i}\mid 0\le i<n \}=r^{-tm}\langle r^d\rangle$.
	\end{itemize}

    Thus,
	$
	\omega_\alpha(G)
	=
	\langle r^d\rangle\cup r^{-tm}\langle r^d\rangle
	\subseteq\langle r\rangle$. We know that $\langle r^d\rangle\cup a\langle r^d\rangle$, for $a\in A$, is a subgroup of $A$ if and only if $a^2\in \langle r^d\rangle$. We have  $n\mid m(k+1)$. Then $d\mid m(k+1)$ . By the definition of $d$, we have $d\mid (1-k)$. Hence $k\equiv 1 \pmod{d}$ and so $k+1\equiv 2\pmod{d}$. Therefore $d\mid 2m$ and so $d\mid 2mt$. This implies that $(r^{-tm})^2\in \langle r^d\rangle$. Thus $\omega_\alpha(G)\le A$.
	\end{proof}

\begin{definition}
	Let $G$ be a group and $H$ a  subgroup of $G$. Then we define $\mathcal{L}(H)$ as follows:
	\[
	\mathcal{L}(H) = 
	\begin{cases} 
		|H|, & \text{if } \omega_\alpha(G)\subseteq H, \\
		\bigl|H\setminus \omega_\alpha(G)\bigr|, & \text{if } \omega_\alpha(G)\nsubseteq H.
	\end{cases}
	\]
\end{definition}

%%%%%%%%%%%%%%%%%%%%%

\begin{lemma}\label{lem:normal}
	Let $G$ be a group and $\alpha\in\Aut(G)$ be an involutory automorphism such that $\omega_\alpha(G)\le G$.
	Let $H$ be a proper normal subgroup of $G$ such that $\alpha(H)=H$. 
	Let $b$ be an integer with
$0\leq b\leq\mathcal{L}(H)$.
Then $H$ is a $(0,b)$-regular set
of $(G,\alpha)$ if and only if
at least one of the following
conditions holds:

\begin{enumerate}
\item $b$ is even;
\item $H$ is a subgroup perfect code
of $(G,\alpha)$.
\end{enumerate}
\end{lemma}

\begin{proof}
Put $\iota(g)=\alpha(g^{-1})$.
Since $H\trianglelefteq G$ and
$\alpha(H)=H$, the involution
$\iota$ permutes the right
cosets of $H$.
For each $x\in G\setminus H$, put
$W_x=Hx\setminus\omega_\alpha(G)$.
By Lemma~\ref{L(x)}, we have
$|W_x|\geq\mathcal L(H)$.
Moreover, $\iota$ permutes
these admissible sets.
Consider the orbits of $\iota$
on the right cosets of $H$ in $G\setminus H$.
If two distinct cosets are
interchanged by $\iota$, choose
any $b$ admissible elements
from one coset and their
images under $\iota$ from
the other.
Now suppose that a right
coset $Hx$ is fixed by $\iota$.
Then $W_x$ is $\iota$-invariant
and decomposes into fixed
points and two-element orbits.
If $b$ is even, we can choose
an $\iota$-invariant subset
of $W_x$ of size $b$.
If $b$ is odd, then $H$
is a subgroup perfect code
by hypothesis.
Since $\iota(Hx)=Hx$, we have
$
\alpha(HxH)=Hx^{-1}H.
$
Also, $H\trianglelefteq G$
implies
$
|HxH|/|H|=1.$
Thus, by Lemma~\ref{ThA},
$Hx\cap\Omega_\alpha(G)
\neq\emptyset.
$
Consequently, $W_x$ contains
an $\iota$-fixed point.
Hence an $\iota$-invariant
subset of $W_x$ of every
cardinality $b\leq |W_x|$
can be selected.
Taking the union of these
selections over all orbits,
we obtain $S_b\subseteq G\setminus H$
such that
$\iota(S_b)=S_b$,
$S_b\cap\omega_\alpha(G)=\emptyset$,
and
$|S_b\cap Hx|=b
\quad\text{for every }x\notin H$.
The result follows from
Lemma~\ref{2.7}.

Conversely, suppose that $H$ is
a $(0,b)$-regular set of
$(G,\alpha)$ for some odd integer
$b$, with connection set $S$.
Assume that $H$ is not a subgroup
perfect code. By Lemma~\ref{ThA},
there exists $x\in G\setminus H$
such that
$\alpha(HxH)=Hx^{-1}H$,
$HxH\cap\Omega_\alpha(G)=\emptyset$,
and $|HxH|/|H|$ is odd.
Since $H\trianglelefteq G$,
we have $HxH=Hx$ and
$\iota(Hx)=Hx$, where
$\iota(g)=\alpha(g^{-1})$.
By Lemma~\ref{2.7},
$|S\cap Hx|=b.
$
Moreover, $S\cap Hx$ is
$\iota$-invariant and contains
no fixed point of $\iota$,
because
$S\cap\omega_\alpha(G)=\emptyset
$
and
$Hx\cap\Omega_\alpha(G)=\emptyset.
$
Thus its elements lie in
two-element orbits, implying
that $b$ is even, a contradiction.
Therefore $H$ is a subgroup
perfect code of $(G,\alpha)$.
\end{proof}

%%%%%%%%%%%
%%%%%%%%%%%
%%%%%%%%% new version %%%%%%%%%%%%%%%%%%%%%%%%%%%

\begin{lemma}\label{regular_b}
Let $G$ be a group having a
cyclic subgroup $A$ of index $2$, and
let $\alpha\in\operatorname{Aut}(G)$ be an involutory automorphism of $G$ such that $\alpha(A)=A$.
Let $H$ be a proper $\alpha$-invariant
subgroup of $G$.
Let $b$ be an integer with
$0\leq b\leq\mathcal{L}(H)$.
Then $H$ is a $(0,b)$-regular set
of $(G,\alpha)$ if and only if
at least one of the following
conditions holds:

\begin{enumerate}
\item $b$ is even;
\item $H$ is a subgroup perfect code
of $(G,\alpha)$.
\end{enumerate}
\end{lemma}

\begin{proof}
By Lemma~\ref{omega_sub_new},
$\omega_\alpha(G)$ is a subgroup of $A$.
Put $\iota(g)=\alpha(g^{-1})$.

If $H\leq A$, then $H\trianglelefteq G$,
and the result follows from
Lemma~\ref{lem:normal}.
Thus we assume that $H\nleq A$.

Put $K=H\cap A$ and choose
$s\in H\setminus A$.
Then $K\trianglelefteq G$,
$H=K\cup Ks$, and $s^2\in K$.
Set $M=K\omega_\alpha(G)$.
Since $K,\omega_\alpha(G)\leq A$, we have
$M\trianglelefteq G$ and
$\alpha(M)=M$. Consequently,
$\iota(M)=M$ and $\iota(sM)=sM$.
Note that every right coset of $H$ is of
the form $Ha$, where $a\in A$, and
$Ha=Ka\;\dot\cup\;Ksa$.
Moreover,
$$
|Ka\setminus \omega_\alpha(G)|=
\begin{cases}
|K|-|K\cap \omega_\alpha(G)|,&a\in M,\\
|K|,&a\notin M,
\end{cases}
$$
and $|Ksa|=|K|$.
Since $K\trianglelefteq G$ and
$\alpha(K)=K$, the involution $\iota$
permutes the right cosets of $K$,
preserving $A$, $G\setminus A$,
$M$, and $sM$.
In particular, if
$\iota(Ka)=Ka'$ or
$\iota(Ksa)=Ksa'$, then
$a\in M$ if and only if $a'\in M$.
Hence the numbers of admissible
elements in the corresponding
$K$-cosets are equal.

$\bullet$ First suppose that $b$ is even.
For each $Ha\ne H$, put
$p=|Ka\setminus \omega_\alpha(G)|$ and $q=|K|$.
By Lemma~\ref{L(x)}, we have
$b\leq \mathcal{L}(H)\leq p+q$.
If $b<p+q$, then
$
b\leq
2\left\lfloor\frac p2\right\rfloor+
2\left\lfloor\frac q2\right\rfloor.
$
Hence there exist even integers
$u,v$ such that
$0\leq u\leq p$,
$0\leq v\leq q$, and $u+v=b$.
We choose $u$ admissible elements
from each $Ka$ and $v$ elements
from each $Ksa$, so that the
resulting sets are preserved by
$\iota$.
This is possible because $\iota$
permutes the admissible $K$-cosets
of each type, preserving their
cardinalities. For two cosets
interchanged by $\iota$, choose
the required elements in one and
their images in the other.
For an $\iota$-invariant coset,
we choose an invariant subset of
the prescribed even cardinality.
If $b=p+q$, simply choose all
admissible elements of $Ha$.
Since the condition $a\in M$
is preserved by $\iota$,
these choices are compatible
throughout $G\setminus H$.
Thus we construct $S_b\subseteq G\setminus H$
such that
$\iota(S_b)=S_b$,
$S_b\cap \omega_\alpha(G)=\emptyset$,
and
$|S_b\cap Ha|=b$, for every $Ha\ne H$.

$\bullet$ Now suppose that $b$ is odd
and $H$ is a subgroup perfect code.
Let $S_1$ be a corresponding
connection set. Then
$\iota(S_1)=S_1$,
$|S_1\cap Ha|=1$, where $Ha\ne H$.
Remove $S_1$ from the admissible
elements. 
Since $S_1$ contains exactly one element
from each $Ha\ne H$, the numbers of
admissible elements in $Ka$ and $Ksa$,
after removing $S_1$, are either
$(p-1,q)$ or $(p,q-1)$, depending on
which part contains the unique element
of $S_1\cap Ha$.
Since $\iota(S_1)=S_1$ and $\iota$
preserves $A$ and $G\setminus A$,
the resulting cardinalities
are preserved whenever the
corresponding $K$-cosets are
interchanged by $\iota$.
Indeed, for every right $K$-coset $R$,
since $\iota(S_1)=S_1$ and
$\iota(\omega_\alpha(G))
=\omega_\alpha(G)$, we have
$
\iota\bigl(
R\setminus(S_1\cup\omega_\alpha(G))
\bigr)
=
\iota(R)\setminus
(S_1\cup\omega_\alpha(G)).
$
Thus the remaining admissible
elements in corresponding
$K$-cosets have the same
cardinality. 
Moreover, since $S_1$ contains exactly
one element from each $Ha$, if
$\iota(Ka)=Ka'$, then the numbers of
remaining admissible elements are equal
not only in $Ka$ and $Ka'$, but also
in $Ksa$ and $Ksa'$.
The same conclusion holds if
$\iota(Ksa)=Ksa'$.
Thus the numbers of remaining
admissible elements in both parts
are constant along any chain of
these correspondences.
Therefore, we can choose the elements
so that the resulting set is preserved
by $\iota$.
As $b-1$ is even and
$b-1\leq \mathcal{L}(H)-1$,
we repeat the construction
for even values, using the
remaining admissible elements.
If $b-1$ equals their total
number in a right coset, we
select all of them; otherwise,
we choose even numbers of
elements from its two parts.
This gives an
$\iota$-invariant set $S'$, where
$S'\cap (S_1\cup \omega_\alpha(G))=\emptyset$,
and
$|S'\cap Ha|=b-1$.
Therefore $S_b=S_1\cup S'$
satisfies
$\iota(S_b)=S_b$,
$S_b\cap \omega_\alpha(G)=\emptyset$,
and
$|S_b\cap Ha|=b$.
The result follows from
Lemma~\ref{2.7}.

Conversely, suppose that $H$ is
a $(0,b)$-regular set of
$(G,\alpha)$ for some odd integer
$b$, with connection set $S$.
Assume that $H$ is not a subgroup
perfect code. By Lemma~\ref{ThA},
there exists $x\in G\setminus H$
such that
$
\alpha(HxH)=Hx^{-1}H,$
$HxH\cap\Omega_\alpha(G)=\emptyset,$
and
$[H:H\cap xHx^{-1}]
$ is odd.
Put $K=H\cap A$.
Since $A$ is cyclic and normal
in $G$, we have $K\trianglelefteq G$
and $[H:K]\leq 2$.
Consequently,
$[H:H\cap xHx^{-1}]\leq 2.
$ As this index is odd, it is equal
to $1$, and hence $x\in N_G(H)$.
Thus $HxH=Hx$ and
$\iota(Hx)=Hx$, where
$\iota(g)=\alpha(g^{-1})$.
By Lemma~\ref{2.7},
$|S\cap Hx|=b.
$ Since $S\cap Hx$ is
$\iota$-invariant,
$Hx\cap\Omega_\alpha(G)=\emptyset$, and
$S\cap\omega_\alpha(G)=\emptyset$, we get that
the involution $\iota$ acts
without fixed points on
$S\cap Hx$.
It follows that $b$ is even,
a contradiction.
Therefore $H$ is a subgroup
perfect code of $(G,\alpha)$.
\end{proof}

%%%%%%%%%%%%%%%%%%%%%%%

\begin{theorem}\label{thm:uni-new}
Let $G$ be a group and let $A$ be a cyclic
subgroup of $G$ of index $2$. Let
$\alpha\in\operatorname{Aut}(G)$ satisfy
$\alpha^2=\mathrm{id}$ and $\alpha(A)=A$.
Let $H$ be a proper subgroup of $G$ such that
$\alpha(H)=H$.
Then, for nonnegative integers $a$ and $b$,
$H$ is an $(a,b)$-regular set of $(G,\alpha)$
if and only if
$
0\le a\le |H\setminus\omega_\alpha(G)|$,
 $0\le b\le\mathcal{L}(H)$,
and the following conditions hold:
\begin{enumerate}
    \item either $a$ is even or
    $H\cap\Omega_\alpha(G)\ne\emptyset$;
    \item either $b$ is even or
    $H$ is a subgroup perfect code of $(G,\alpha)$.
\end{enumerate}
\end{theorem}

    \begin{proof}
By Lemma~\ref{omega_sub_new},
$\omega_\alpha(G)$ is a subgroup
of $A$.
The result follows immediately
from Lemmas~\ref{a-regularity},
\ref{regular_b}, and~\ref{2.7}.
\end{proof}
%%%%%%%%%%%%%%%

\begin{corollary}\label{cor:characteristic}
Let $G$ be a group with a
characteristic cyclic subgroup of index $2$.
Then, for every $\alpha\in\operatorname{Aut}(G)$
with $\alpha^2=\mathrm{id}$ and every proper
$\alpha$-invariant subgroup $H$ of $G$,
the results of Theorem~\ref{thm:uni-new}
hold.
\end{corollary}

\begin{proof}
Let $A$ be a characteristic cyclic subgroup
of $G$ of index $2$.
Then
$\alpha(A)=A$ for every
$\alpha\in\operatorname{Aut}(G)$.
Hence the result follows from
Theorem~\ref{thm:uni-new}.
\end{proof}

In particular, the previous corollary applies to several classes of groups, including dihedral groups $D_{2n}$ ($n\ge 3$), dicyclic groups $\mathrm{Dic}_m$ ($m\ge 3$), and semidihedral groups $SD_{2^n}$ ($n\ge 4$), among others. Moreover, Theorem~\ref{thm:uni-new} also applies to modular groups $M_{2^n}$ ($n\ge 4$), provided that their cyclic subgroup of index $2$ is $\alpha$-invariant.
This naturally leads us to ask:

\begin{problem}
 For which finite groups $G$ is $\omega_\alpha(G)$ a subgroup of $G$ for every involutory automorphism $\alpha$?
 \end{problem}

\begin{problem}
If $\omega_\alpha(G)$ is a subgroup
of $G$, under what additional
conditions does the characterization
of $(a,b)$-regular subgroups in
Theorem~\ref{thm:uni-new}
remain valid?
\end{problem}

\end{document}